\documentclass[11pt,showlabels]{article1}
\usepackage{amssymb}
\usepackage{epsfig}
\newtheorem{theorem}{Theorem}
\newtheorem{lemma}{Lemma}
\newtheorem{remark}{Remark}
\newtheorem{proposition}{Proposition}

\title{\bf Optimal Control for Discrete-time Markov Jump Linear System with Control Input Delay\thanks{This work is supported by the National Natural Science
Foundation of China (Nos. 61473134, 61573220, 61120106011, 61573221) and the Postdoctoral Science Foundation of China (No.
2017M622231). $^{*}$Corresponding author: Huanshui Zhang. Email: hszhang@sdu.edu.cn}
}

\author{Chunyan Han$^a$,\ Hongdan Li$^b$, \ Huanshui Zhang$^{b,*}$
\ \\
\\
\ \ \ $^a$ School of Electrical Engineering, University of Jinan, \\Jinan Shandong 250022, China
\\
$^b$ School of Control Science and Engineering, Shandong
University, \\Jinan Shandong 250061, China}

\begin{document}
\baselineskip 16pt
\date{}
  \maketitle
\begin{abstract}
This paper deals with the finite horizon optimal control problem for discrete-time Markov jump linear system with input delay. The correlation among the jumping parameters and the input delay are considered simultaneously, which forms the basic difficulty of the design. one of the key techniques is to solve a delayed forward and backward jumping parameter difference equation which is obtained by an improved maximum principle, and the other is the introduction of a ``d-step backward formula". Based on the proposed techniques, a necessary and sufficient condition for the existence of the optimal controller is given in an explicit form and an analytical solution to the optimal controller is supplied. The optimal controller is a linear function of the current time state and the historical time control input, where the feedback gains are a set of jumping parameter matrices derived by solving a new type of
 coupled difference Riccati equation. The key step in the derivation is to establish the relationship between the costate and the real state of the system. The result obtained in this paper can be viewed as a generalization of the standard case, in which there is only one mode of operation.

%
\bigskip

\end{abstract}

\pagestyle{plain} \setcounter{page}{1}
\section{Introduction}

Optimal control and stabilization of Markov jump linear systems (MJLSs) have matured over the last decades \cite{BS75}-\cite{CF16}. Such systems are encountered in many areas of control engineering including aircraft and nuclear power control, flight control, and etc \cite{Mar1990}. In the previous research works, Costa has made great contributions to the study of this problem and many celebrated results have been proposed in \cite{CABM99}-\cite{CF16}. The fundamental tool for dealing with optimal control in presence of jumping parameters is the dynamic programming method. However, when time delays in the control input are present, this method can't be applied directly since the historical input terms are existed in the step-by-step design.

Optimal control of dynamic systems with input delays has received growing attention in recent years due to novel application areas such as control over networks \cite{LYZ11}-\cite{YLT14}. For the deterministic system with input delay, the reader may refer to the excellent surveys in \cite{WI81}-\cite{ZDX06}. In \cite{WI81}, the control problem for a single deterministic input-delayed system was considered, and an optimal controller was developed by the Smith predictor. In \cite{Art82}, a reduction technique was developed to convert the stabilization problem subject to input delays to an equivalent delay-free one. To overcome the limitation of the original Smith predictor and the model reduction technique, some renewed prediction approaches, such as truncated predictor feedback \cite{ZLD12} and closed-loop predictor approaches\cite{ZLM17}, \cite{CG17}, have been developed. In \cite{ZDX06}, an alternative and efficient approach-the duality method, has been developed for the studying on linear quadratic regulation of systems with multiple input delays. Obviously, the aforementioned works have supplied good results for the advances of optimal control theories on deterministic delayed systems. However, the results can't be applied to the delayed MJLS directly since the separation principle is satisfied for deterministic case but not suitable for stochastic case.

Not requiring the satisfaction of separation principle, in \cite{ZWL12}-\cite{LZ16}, a direct approach based on a  new stochastic maximum principle was proposed for dealing with the optimal control of system with multiplicative noise and input delay simultaneously. A necessary and sufficient condition for the existence of the optimal controller was supplied, and an explicit solution was given. These results paved new ways for the investigation of optimal control for stochastic systems with input delays. Motivated by this, we consider the optimal control problems for MJLS with input delay in this paper. Compared to the multiplicative noise system, the MJLS is more complicated due to the correlation of jumping parameters. And thus some improved and additional techniques need to be developed further. To the best of our knowledge, the study of optimal control for MJLS with input delay hasn't been reported before.

In this paper, motivated by \cite{ZWL12}-\cite{LZ16}, we investigate the optimal control for MJLS with input delay. The simultaneous appearance of jumping parameters and input delay forms the fundamental difficulty in investigation. Two basic formulas are developed: one is an improved delayed forward and backward jumping parameter equation (D-FBJPE) which is used to deal with the input delay, and the other is a d-step backward formula which is used to overcome the correlation of the jumping parameters. Based on the the two proposed techniques, a necessary and sufficient condition for the existence of the optimal controller is developed, and an explicit solution is given for the first time in terms of a new type of coupled difference Riccati equation. The key step in the derivation is to establish the relationship between the optimal state and the costate.

In section
2, we present the problem formulations. In
section 3, some preparations are given and the main results are proposed. A numerical example is given in section 4 and some conclusion remarks are made in section 5.

Notations: Throughout this paper, ${R}^n$ denotes the $n$-dimensional Euclidean space,
$R^{m\times n}$ denotes the norm bounded linear space of all
$m\times n$ matrices.
For $L\in R^{n\times n}$, $L'$ stands for the transpose of $L$. As usual, $L\geq 0 (L>0)$
will mean that the symmetric matrix $L\in R^{n\times n}$ is positive semi-definite
(positive definite), respectively.
 $\mbox{E}(.)$ denotes the mathematical expectation operator, $\mbox{P}(.)$ means the occurrence probability of an event.
We will compactly write the sum $\sum_{l_{k-d+1}=1}^L\lambda_{l_{k-d},l_{k-d+1}}~\cdots ~\sum_{l_{k}=1}^L\lambda_{l_{k-1},l_k}$ as $\Lambda_{l_{k-d},l_{k}}$.


\section{Problem Statement}
\setcounter{equation}{0}
We consider in this paper the finite horizon
optimal control for MJLS with input delay when the state variable $x(k)$ and the jump variable $\theta(k)$ are available to the controller. On the stochastic basis $(\Omega,{\cal{G}}, {\cal{G}}_k,\mbox{P})$, consider the following MJLS with input delay
\begin{eqnarray}
x(k+1)&=&A_{\theta(k)}(k)x(k)+B_{\theta(k)}(k)u(k-d),\label{fff1}
\end{eqnarray}
where $x(k)\in {\mbox{R}}^n$ is the state, $u(k)\in {\mbox{R}}^m$  is the control input with delay $d>0$. $\theta(k)$ is a discrete-time Markov chain with finite state space $\Theta\triangleq\{1,2,\cdots,L\}$ and transition probability $\lambda_{l_k,l_{k+1}}=\mbox{P}(\theta(k+1)=l_{k+1}|\theta(k)=l_k)(l_k,l_{k+1}=1,2,\cdots,L)$. We set $\pi_{l_k}(k)=\mbox{P}(\theta(k)=l_k)(l_k=1,2,\cdots,L)$, while $A_{l_k}(k), B_{l_k}(k)({l_k}=1,\cdots,L)$ are matrices of appropriate dimensions. The initial values $x_0,u(i),i=-d,\cdots,-1$ are known. We assume that $\theta(k)$ is independent of $x_0$ and $u(i),i=-d,\cdots,-1$.

The quadratic cost
associated to system (\ref{fff1}) with admissible control law $u=(u(0),\cdots,u(N-d))$ is given by
\begin{eqnarray}
J_N&=&\mbox{E}[\sum_{k=0}^Nx(k)'Q_{\theta(k)}(k)x(k)+\sum_{k=d}^Nu(k-d)'R_{\theta(k)}(k)u(k-d)\nonumber\\
&&+x(N+1)'P_{\theta(N+1)}(N+1)x(N+1)],\label{f2}
\end{eqnarray}
where $N>d$ is an integer, $x(N+1)$ is the terminal state, $P_{l_{N+1}}(N+1)(l_{N+1}=1,\cdots,L)$ reflects the penalty on the terminal state, the matrix functions $R_{l_k}(k)\geq0(l_k=1,\cdots,L)$ and $Q_{l_k}(k)\geq 0(l_k=1,\cdots,L)$. Denote ${\cal{G}}_k=\{\theta(t);t=0,\cdots,k\}$, so the problem considered in this paper can be stated as:

\textbf{Problem 1}: Find a ${\cal{G}}_k$-measurable controller $u(k)$ such that (\ref{f2}) is minimized subject to (\ref{fff1}).

\begin{remark}
 For brevity, we will omit the time steps in the system matrices and the penalty matrices in the following discussions. That is denoting $A_{\theta(k)}(k), B_{\theta(k)}(k)$, $Q_{\theta(k)}(k), R_{\theta(k)}(k)$ and $P_{\theta(N+1)}(N+1)$ as $A_{\theta(k)}, B_{\theta(k)}$, $Q_{\theta(k)}, R_{\theta(k)}$ and $P_{\theta(N+1)}$, respectively. This will not affect the final results.
 \end{remark}
\begin{remark}
It should be pointed out that the result in \cite{ZLXF15} can not be applied to the optimal control for MJLS with $d$-step ($d\geq 2$) input delay directly since the dependence of the jumping parameters. A new version of maximum principle needs to be developed and a ``d-step backward formula" will be employed which form the basic tools for the design of the optimal controller.
\end{remark}

\section{Optimal Control for MJLS with input delay}
In this section, a necessary and sufficient condition for the existence of the optimal controller is given, and the explicit expressions for the optimal controller, the optimal costate, and the optimal cost are derived. The key techniques employed in this part are the Markovian jump maximum principle and the ``d-step backward formula".
\subsection{Preparations}
Due to the dependence of $\theta(k)$ on its past values, the new version of the maximum principle for the optimal control with jumping parameters and input dleay needs to be established.
\begin{lemma} (maximum principle)
According to system (\ref{fff1}) and performance (\ref{f2}). If problem 1 is solvable, then the optimal ${\cal{G}}_{k-d}$-measurable controller $u(k-d)$ satisfies
\begin{eqnarray}
0=\mbox{E}[B_{\theta(k)}'\lambda_k+R_{\theta(k)}u(k-d)|{\cal{G}}_{k-d}],k=d,\cdots,N,\label{f3}
\end{eqnarray}
where the costate $\lambda_k$ satisfies
\begin{eqnarray}
\lambda_N&=&\mbox{E}[P_{\theta(N+1)}x(N+1)|{\cal{G}}_N],\label{f4}\\
\lambda_{k-1}&=&\mbox{E}[A_{\theta(k)}'\lambda_k+Q_{\theta(k)}x(k)|{\cal{G}}_{k-1}],k=0,\cdots,N.\label{f5}
\end{eqnarray}
\end{lemma}

\emph{Proof}. See Appendix A.

From (\ref{fff1}), we get that
\begin{eqnarray}
x(k)=F_{\theta(k-1),\theta(k-d)}x(k-d)+\sum_{i=k-d}^{k-1}F_{\theta(k-1),\theta(i+1)}B_{\theta(i)}u(i-d),\label{fff2}
\end{eqnarray}
where $F_{\theta(k-1),\theta(k-d)}, F_{\theta(k-1),\theta(i+1)}$ are as in (\ref{f8}). Based on the transformation (\ref{fff2}), we get the following expression.
\begin{lemma} (d-step backward formula)
For any ${\cal{G}}_k-$measurable function $f_{\theta(k)}$, we have the following relation
\begin{eqnarray}
\mbox{E}\{f_{\theta(k)}x(k)|{\cal{G}}_{k-d}\}=\Lambda_{l_{k-d},l_k}[f_{l_k}F_{l_{k-1},l_{k-d}}]x(k-d)+\sum_{i=k-d}^{k-1}\Lambda_{l_{k-d},l_{k}}[f_{l_k}F_{l_{k-1},l_{i+1}}B_{l_i}]u(i-d).\label{fff3}
\end{eqnarray}
\end{lemma}
\emph{Proof}. The result can be obtained directly and thus the proof is omitted here.
\begin{remark}
In Lemma 2, $x(k)$ moved $d$ steps back, so we named (\ref{fff3}) as a $d$-step backward formula. However, it should be pointed out that
\begin{eqnarray}
\mbox{E}\{f_{\theta(k)}x(k)|{\cal{G}}_{k-d}\}\neq \mbox{E}\{f_{\theta(k)}|{\cal{G}}_{k-d}\}\mbox{E}\{x(k)|{\cal{G}}_{k-d}\}\label{fff4}
\end{eqnarray}
since $x(k)=A_{\theta(k-1)}x(k-1)+B_{\theta(k-1)}u(k-1-d)$ and $\theta(k)$ is correlated with $\theta(k-1)$. If system (\ref{fff1}) becomes the multiplicative noise system investigated in \cite{ZLXF15}
\begin{eqnarray}
x(k+1)=(A+\omega_k\bar{A})x(k)+(B+\omega_k\bar{B})u(k-d)\nonumber
\end{eqnarray}
and $f_{\theta(k)}=\omega_k f$ with $\omega_k$ being a multiplicative noise, then the equality
\begin{eqnarray}
\mbox{E}\{f_{\theta(k)}x(k)|{\cal{G}}_{k-d}\}= \mbox{E}\{f_{\theta(k)}|{\cal{G}}_{k-d}\}\mbox{E}\{x(k)|{\cal{G}}_{k-d}\}.\nonumber
\end{eqnarray}
holds \cite{ZLXF15}. It can be seen that the optimal control for MJLS with input delay is more complicated than that of multiplicative noise system with input delay, and (\ref{fff4}) is the key difficulty existed in the derivation. Fortunately, we reveal the relationship (\ref{fff3}), which forms the basic formula for the controller design for the delayed MJLS.
\end{remark}


\subsection{Solution to the optimal control}

In this part, we will derive the analytic solution to the optimal control for MJLS with input delay by Lemma 1 and Lemma 2.

Firstly, we introduce some notations $W_{l_{k-d}}(k-d), T_{l_{k-d}}^j(k-d) (j=0, 1, \cdots, d-1)$, $P_{l_{k-1}}(k-1), P_{l_{k-1}}^0(k-1), \delta_{l_{k-1}}^j(k-1)(j=1, 2, \cdots, d-1)$ and $\alpha_{l_{k-1},l_{k-j}}^{d-j}(k-1,k-j)(j=1, 2, \cdots, d-1)$ for $k=N, N-1,\cdots,0$, $\theta(k-d)=l_{k-d}\in \Theta,\theta(k-1)=l_{k-1}\in \Theta$. To save space, we rewrite the preceding notations as $W_{l_{k-d}}, T_{l_{k-d}}^j(j=0, 1, \cdots, d-1)$, $P_{l_{k-1}}, P_{l_{k-1}}^0, \delta_{l_{k-1}}^j(j=1, 2, \cdots, d-1)$ and $\alpha_{l_{k-1},l_{k-j}}^{d-j}(j=1, 2, \cdots, d-1)$ respectively.

In this part, $W_{l_{k-d}}$ and $T_{l_{k-d}}^j(j=0, 1, \cdots, d-1)$ satisfy the following backward equations
\begin{eqnarray}
W_{l_{k-d}}&=&\Lambda_{l_{k-d},l_k}[B_{l_k}'(P_{l_k}-P_{l_k}^0)B_{l_k}+R_{l_k}] -\sum_{s=1}^{d-1}\{\Lambda_{l_{k-d},l_{k-s}}[(T_{l_{k-s}}^s)'W_{l_{k-s}}^{-1}
T_{l_{k-s}}^s]\},\label{a1}\\
T_{l_{k-d}}^0&=&\Lambda_{l_{k-d},l_k}[B_{l_k}'(P_{l_k}-P_{l_k}^0)F_{l_k,l_{k-d+1}}]-\sum_{s=1}^{d-1}\{\Lambda_{l_{k-d},l_{k-s}}[(T_{l_{k-s}}^s)'W_{l_{k-s}}^{-1}T_{l_{k-s}}^0\nonumber\\
&&\times F_{l_{k-s},l_{k-d+1}}]\},\label{a2}\\
T_{l_{k-d}}^1&=&\Lambda_{l_{k-d},l_k}[B_{l_k}'(P_{l_k}-P_{l_k}^0)F_{l_{k},l_{k-d+2}}B_{l_{k-d+1}}]-\sum_{s=1}^{d-2}\{\Lambda_{l_{k-d},l_{k-s}}[(T_{l_{k-s}}^s)'W_{l_{k-s}}^{-1}T_{l_{k-s}}^0\nonumber\\
&&\times F_{l_{k-s},l_{k-d+2}}B_{l_{k-d+1}}]\}-\Lambda_{l_{k-d},l_{k-d+1}}[(T_{l_{k-d+1}}^{d-1})'W_{l_{k-d+1}}^{-1}T_{l_{k-d+1}}^0B_{l_{k-d+1}}],\label{a3}\\
T_{l_{k-d}}^j&=&\Lambda_{l_{k-d},l_k}[B_{l_k}'(P_{l_k}-P_{l_k}^0)F_{l_k,l_{k-d+j+1}}B_{l_{k-d+j}}]-\sum_{s=1}^{d-j}\{\Lambda_{l_{k-d},l_{k-s}}[(T_{l_{k-s}}^s)'W_{l_{k-s}}^{-1}\nonumber\\
&&\times T_{l_{k-s}}^0F_{l_{k-s},l_{k-d+j+1}}B_{l_{k-d+j}}]\}-\sum_{s=d-j+1}^{d-1}\{\Lambda_{l_{k-d},l_{k-s}}[(T_{l_{k-s}}^s)'W_{l_{k-s}}^{-1} T_{l_{k-s}}^{s-(d-j)}]\}, \nonumber\\
&&j=2,3,\cdots,d-1 \label{a4}
\end{eqnarray}
for $k=N, N-1,\cdots,0, l_{k-d}\in \Theta$ with terminal values
\begin{eqnarray}
T_{l_{N-i}}^j=0, j=0,1,\cdots,d-1, i=0,1,\cdots,d-1, l_{N-i}\in \Theta\nonumber
\end{eqnarray}
and $P_{l_{k-1}}$ and $P_{l_{k-1}}^0$ satisfy the following backward recursions
\begin{eqnarray}
P_{l_{k-1}}&=&\Lambda_{l_{k-1},l_k}[Q_{l_k}+A_{l_k}'(P_{l_k}-P_{l_k}^0)A_{l_k}],\label{a5}\\
P_{l_{k-1}}^0&=&(T_{l_{k-1}}^0)'W_{l_{k-1}}^{-1}T_{l_{k-1}}^0\label{a6}
\end{eqnarray}
for $k=N, N-1,\cdots,0, l_{k-1}\in \Theta$ with terminal values
\begin{eqnarray}
P_{l_N}&=&\Lambda_{l_{N},l_{N+1}}P_{l_{N+1}},\nonumber\\
P_{l_{N-i}}^0&=&0, i=0,1,\cdots,d-1, l_{N-i}\in \Theta.\nonumber
\end{eqnarray}
(\ref{a1})-(\ref{a6}) is termed as the backward coupled difference Riccati equation, which is with the same dimension as that of the original system state. An existence condition and an explicit solution to the optimal controller will be given in terms of the Riccati equation.

In addition, the introduction of $\delta_{l_{k-1}}^j(j=1, 2, \cdots, d-1)$ and $\alpha_{l_{k-1},l_{k-j}}^{d-j}(j=1, 2, \cdots, d-1)$ is critical to establish the relationship between the optimal original state $x(k)$ and costate $\lambda_{k-1}$, where $\delta_{l_{k-1}}^j(j=1, 2, \cdots, d-1)$ and $\alpha_{l_{k-1},l_{k-j}}^{d-j}(j=1, 2, \cdots, d-1)$ satisfy the following expressions
\begin{eqnarray}
(\delta_{l_{k-1}}^1)'&=&\Lambda_{l_{k-1},l_k}[A_{l_k}'(P_{l_k}-P_{l_k}^0)B_{l_k}]-(T_{l_{k-1}}^0)'W_{l_{k-1}}^{-1}T_{l_{k-1}}^1,\label{a7}\\
(\delta_{l_{k-1}}^j)'&=&\Lambda_{l_{k-1},l_k}[A_{l_k}'(\delta_{l_k}^{j-1})']-(T_{l_{k-1}}^0)'W_{l_{k-1}}^{-1}T_{l_{k-1}}^j,\nonumber\\
&& j=2,3,\cdots,d-1, \label{a8}\\
(\alpha_{l_{k-1},l_{k-1}}^{d-1})'&=&(\delta_{l_{k-1}}^{d-1})',\label{a9}\\
(\alpha_{l_{k-1},l_{k-j}}^{d-j})'&=&(\delta_{l_{k-1}}^{d-j})'-\sum_{s=1}^{j-1}(\alpha_{l_{k-1},l_{k-s}}^{d-s})'W_{l_{k-s-1}}^{-1}T_{l_{k-s-1}}^{d-j+s},\label{a10}\\
&&j=2,3,\cdots,d-1,\nonumber
\end{eqnarray}
for $k=N, N-1,\cdots,0, l_{k-1}, l_{k-j}\in \Theta$. Moreover, The relationship between $T_{l_{k-d}}^j(j=0, 1, \cdots, d-1)$ and $\alpha_{l_{k-1},l_{k-j}}^{d-j}(j=1, 2, \cdots, d-1)$ is
established in the following proposition, which will be used in the derivation of the main result.
\begin{proposition}
Consider $\alpha_{l_{k-1},l_{k-j}}^{d-j}$ and $T_{l_{k-j}}^j$ as in (\ref{a2})-(\ref{a4}), (\ref{a9}) and (\ref{a10}), the following expressions are satisfied
\begin{eqnarray}
\mbox{E}\{A_{l_k}'(\alpha_{l_k,l_k}^{d-1})'|{\cal{G}}_{k-1}\}&=&(T_{l_{k-1}}^0)'\label{c1}\\
\mbox{E}\{A_{l_k}'(\alpha_{l_k,l_{k-j+1}}^{d-j})'|{\cal{G}}_{k-1}\}&=&(\alpha_{l_{k-1},l_{k-j+1}}^{d-j+1})', j=2,\cdots,d-1,\label{c2}\\
\mbox{E}\{B_{l_k}'(\alpha_{l_k,l_k}^{d-1})'|{\cal{G}}_{k-1}\}&=&(T_{l_{k-1}}^1)'\label{c3}\\
\mbox{E}\{B_{l_k}'(\alpha_{l_k,l_{k-j+1}}^{d-j})'|{\cal{G}}_{k-j}\}&=&(T_{l_{k-j}}^j)', j=2,\cdots,d-1.\label{c4}
\end{eqnarray}
\end{proposition}
\emph{Proof}. See Appendix B.
\begin{theorem}
Problem 1 admits a unique optimal control if and only if
\begin{eqnarray}
W_{l_{k-d}}>0, \ k=N,N-1,\cdots,d, l_{k-d}\in \Theta.\label{a12}
\end{eqnarray}
In this case, the analytical solution to the optimal control is given by
\begin{eqnarray}
u(k-d)&=&-W_{l_{k-d}}^{-1}T_{l_{k-d}}^0x(k-d+1)-\sum_{j=1}^{d-1}W_{l_{k-d}}^{-1}T_{l_{k-d}}^ju(k-2d+j)\label{a13}\\
&&k=d,d+1,\cdots,N, l_{k-d}\in \Theta.\nonumber
\end{eqnarray}
The optimal costate is
\begin{eqnarray}
\lambda_{k-1}&=&(P_{l_{k-1}}-P_{l_{k-1}}^0)x(k)-\sum_{s=1}^{d-1}(\alpha_{l_{k-1},l_{k-s}}^{d-s})'W_{l_{k-s-1}}^{-1}\mbox{E}\{\alpha_{l_{k-1},l_{k-s}}^{d-s}x(k)|{\cal{G}}_{k-s-1}\}.\label{a14}
\end{eqnarray}
and the optimal cost is
\begin{eqnarray}
J_N^{*}&=&\mbox{E}\{\sum_{k=0}^{d-1}x(k)'Q_{l_k}x(k)+x(d)'(P_{l_{d-1}}-P_{l_{d-1}}^0)x(d)\nonumber\\
&&-x(d)'\sum_{s=1}^{d-1}(\alpha_{l_{d-1},l_{d-s}}^{d-s})'W_{l_{d-s-1}}^{-1}\mbox{E}[\alpha_{l_{d-1},l_{d-s}}^{d-s}x(d)|{\cal{G}}_{d-s-1}]\}.\label{ff24}
\end{eqnarray}

\end{theorem}

\emph{Proof}. See Appendix C.

\begin{remark}
The solution to the optimal controller is based on a set of a generalized coupled Riccati equations. It can be found that if there is no time delays in the input, the coupled Riccati equation will degenerate into the following one
\begin{eqnarray}
\Upsilon_{l_{k}}&=&[B_{l_k}'(\Lambda_{{l_k},{l_{k+1}}}P_{l_{k+1}})B_{l_k}+R_{l_k}],\nonumber\\
M_{l_{k}}&=&[B_{l_k}'(\Lambda_{{l_k},{l_{k+1}}}P_{l_{k+1}})A_{l_k}],\nonumber\\
P_{l_{k}}&=&A_{l_k}'(\Lambda_{{l_k},{l_{k+1}}}P_{l_{k+1}})A_{l_k}+Q_{l_k}-M_{l_{k}}'\Upsilon_{l_{k}}^{-1}M_{l_{k}},\nonumber\\
&&l_k=1,\cdots,L,\nonumber
\end{eqnarray}
which has been developed in \cite{CFM05}.
\end{remark}

\section{Numerical Examples}
In this part, we present a simple example to illustrate the
theoretical result for the optimal control of MJLS with input delay. Consider a second-order dynamic system (\ref{fff1}) with the performance ({\ref{f2}}). The specifications of the system and the weighting matrices are as follows
\begin{eqnarray}
&&A_1=\left[
      \begin{array}{cc}
        2 & 1.1 \\
        -1.7 & -0.8 \\
      \end{array}
    \right], A_2=\left[
                    \begin{array}{cc}
                      0.8 & 0 \\
                      0 & 0.6 \\
                    \end{array}
                  \right], B_1=\left[
                                  \begin{array}{c}
                                    1 \\
                                    1 \\
                                  \end{array}
                                \right],B_2=\left[
                                               \begin{array}{c}
                                                 2 \\
                                                 1 \\
                                               \end{array}
                                             \right],\nonumber\\
&&Q_{1}=\left[
          \begin{array}{cc}
            1& 0 \\
            0 & 1 \\
          \end{array}
        \right],Q_{2}=\left[
                          \begin{array}{cc}
                            1 & 0 \\
                            0 & 1 \\
                          \end{array}
                        \right],R_{1}=1,R_2=2,P_{N+1}=\left[
          \begin{array}{cc}
            1& 0 \\
            0 & 1 \\
          \end{array}
        \right].\nonumber
\end{eqnarray}
$\theta(k)$ is the Markov chain taking values in a finite set $\{1,\ 2\}$
with transition probability matrix $\left[
                                      \begin{array}{cc}
                                        0.9 & 0.1 \\
                                        0.3 & 0.7 \\
                                      \end{array}
                                    \right]$. The initial distribution of $\theta(k)$ is $(0.5,\ 0.5)$.
The initial
value $x(0)=[2\ 2]', u(-1)=-1,u(-2)=-2$, and the time delay $d=2$.

In this example, the time horizon is set to $N=7$. Without loss of generality, we run $50$
Monte Carlo simulations from $k=0$ to $7$, and select the first trajectory to show the efficiency of the proposed algorithm. By applying Theorem 1, the calculation result for $W_i(k), T^0_i(k), T^1_i(k)(i=1,2)$ are listed in Table 1.
\begin{table}
\centering
\caption{Calculation results}\vspace{3mm}
\begin{tabular}{c|c|c|c|c|c|c}
  \hline
  $k$ & $W_1(k)$ & $W_2(k)$ & $T_1^0(k)$ & $T_2^0(k)$ & $T_1^1(k)$ & $T_2^1(k)$ \\
  \hline
  0 & 23.6031 & 26.7636 & [12.2690\ 7.5948] & [9.6518\ 4.6516] & 21.8683 & 24.7279 \vspace{0.5mm}\\
  1 & 23.1641 & 26.2088 & [12.0539\ 7.4614] & [9.4635\ 4.5596] & 21.4732 & 24.2257 \vspace{0.5mm}\\
  2 & 21.8477 & 24.0482 & [11.6367\ 7.1986] & [8.9148\ 4.2748] & 20.5775 & 22.5743 \vspace{0.5mm}\\
  3 & 17.7981 & 19.0574 & [9.6188\ 5.9405] & [7.1852\ 3.4079] & 16.8338 & 17.9382 \vspace{0.5mm}\\
  4 & 3.6400 & 5.0800 & [0.3659\ 0.2187] & [0.7673\ 0.2769] & 0.9770 & 2.2790\vspace{0.5mm} \\
  5 & 1.1000 & 1.7000 & [0\ 0] & [0\ 0] & 0 & 0 \\
  \hline
\end{tabular}
\end{table}

Note that $W_i(k)>0$ for $i=1,2,k=1,\cdots,5$, thus there exist a unique solution to Problem 1. The optimal controller is computed by
\begin{eqnarray}
u(0)&=&-\left[
          \begin{array}{cc}
            0.3606 & 0.1738 \\
          \end{array}
        \right]x(1)-0.9239u(-1),\nonumber\\
u(1)&=&-\left[
          \begin{array}{cc}
            0.3611 & 0.1740 \\
          \end{array}
        \right]x(2)-0.9243u(0),\nonumber\\
u(2)&=&-\left[
          \begin{array}{cc}
            0.5326 & 0.3295 \\
          \end{array}
        \right]x(3)-0.9419u(1),\nonumber\\
u(3)&=&-\left[
          \begin{array}{cc}
            0.5404 & 0.3338 \\
          \end{array}
        \right]x(4)-0.9458u(2),\nonumber\\
u(4)&=&-\left[
          \begin{array}{cc}
            0.1005 & 0.0601 \\
          \end{array}
        \right]x(5)-0.2684u(3),\nonumber\\
u(5)&=&0,\nonumber
\end{eqnarray}
and the optimal value of (\ref{f2}) is $J_N^*=93.7285$.
%

%
\section{Conclusion}
This paper dealt with the optimal control for MJLS with multi-step input delay. A necessary and sufficient condition for the existence of a unique solution has been developed and a dynamic Markovian jump controller has been given in terms of a coupled difference Riccati equation. One of the key techniques employed in this paper is the maximum principle, the other is the $d$-step backward formula. It should be noted that our derivations avoid the augmented argument, mainly establish and take advantage of the link between the optimal state and the auxiliary variable. Compared with the result for the system with multiplicative noise and input delay, the optimal control for MJLS with input delay is more complicated due to the correlation of the jumping parameters. This is the reason why the problem has not been solved before. The stabilization for MJLS with input delay will be studied in a forthcoming paper.

\appendix
\section{Proof of Lemma 1}
\emph{Proof}. Denote $N$ as the final time. Consider the increment of the control variable $u(k-d)$ and deduce an expression of the corresponding variation of (\ref{f2})
\begin{eqnarray}
dJ_N&=&\mbox{E}[2x(N+1)'P_{\theta(N+1)}dx(N+1)+2\sum_{k=0}^Nx(k)'Q_{\theta(k)}dx(k)\nonumber\\
&&+2\sum_{k=d}^Nu(k-d)'R_{\theta(k)}du(k-d)].\label{f6}
\end{eqnarray}
In view of system (\ref{fff1}), we have
\begin{eqnarray}
dx(k+1)=F_{\theta(k),\theta(0)}dx_0+\sum_{i=0}^kF_{\theta(k),\theta(i+1)}B_{\theta(i)}du(i-d),\label{f7}
\end{eqnarray}
where
\begin{eqnarray}
F_{\theta(k),\theta(i)}&=&A_{\theta(k)}\cdots A_{\theta(i)},i=0,\cdots,k,\nonumber\\
F_{\theta(k),\theta(k+1)}&=&I.\label{f8}
\end{eqnarray}
Plugging the equation (\ref{f7}) in (\ref{f6}) we deduce that
\begin{eqnarray}
dJ_N
&=&\mbox{E}\{2x(N+1)'P_{\theta(N+1)}[F_{\theta(N),\theta(0)}dx_0+\sum_{i=0}^NF_{\theta(N),\theta(i+1)}B_{\theta(i)}du(i-d)]\nonumber\\
&&+2\sum_{i=d}^Nu(i-d)'R_{\theta(i)}du(i-d)+2\sum_{k=0}^Nx(k)'Q_{\theta(k)}F_{\theta(k-1),\theta(0)}dx_0\nonumber\\
&&+2\sum_{i=0}^{N-1}\sum_{k={i+1}}^Nx(k)'Q_{\theta(k)}F_{\theta(k-1),\theta(i+1)}B_{\theta(i)}du(i-d)\}.\label{f11}
\end{eqnarray}
Since we just pay attention to the increment of $J_N$ caused by the increment of $u(i-d)$, the initial state $x_0$ is fixed and its increment $dx_0$ is thus $0$. Therefore,
\begin{eqnarray}
dJ_N
&=&\mbox{E}\{2[x(N+1)'P_{\theta(N+1)}F_{\theta(N),\theta(N+1)}B_{\theta(N)}+u(N-d)'R_{\theta(N)}]du(N-d)\nonumber\\
&&+2\sum_{i=0}^{N-1}[x(N+1)'P_{\theta(N+1)}F_{\theta(N),\theta(i+1)}B_{\theta(i)}+u(i-d)'R_{\theta(i)}\nonumber\\
&&+\sum_{k=i+1}^Nx(k)'Q_{\theta(k)}F_{\theta(k-1),\theta(i+1)}B_{\theta(i)}]du(i-d)\}.\label{f12}
\end{eqnarray}
Define
\begin{eqnarray}
\lambda_i=\mbox{E}\{\sum_{k=i+1}^NF_{\theta(k-1),\theta(i+1)}'Q_{\theta(k)}x(k)+F_{\theta(N),\theta(i+1)}'P_{\theta(N+1)}x(N+1)|{\cal{G}}_i\},\label{f13}
\end{eqnarray}
then we have
\begin{eqnarray}
\lambda_{i-1}
&=&\mbox{E}\{Q_{\theta(i)}x(i)+A_{\theta(i)}'\lambda_i|{\cal{G}}_{i-1}\}\nonumber.
\end{eqnarray}
It has been shown (\ref{f4}) and (\ref{f5}).
Based on (\ref{f13}), we deduce that
\begin{eqnarray}
dJ_N&=&\mbox{E}\{2\sum_{i=d}^N\mbox{E}[u(i-d)'R_{\theta(i)}du(i-d)|{\cal{G}}_i]\nonumber\\
&&+2\mbox{E}[x(N+1)'P_{\theta(N+1)}F_{\theta(N),\theta(N+1)}|{\cal{G}}_N]B_{\theta(N)}du(N-d)\nonumber\\
&&+2\sum_{i=0}^{N_1}\mbox{E}[x(N+1)P_{\theta(N+1)}F_{\theta(N),\theta(i+1)}\nonumber\\
&&+\sum_{k=i+1}^Nx(k)'Q_{\theta(k)}F_{\theta(k-1),\theta(i+1)}|{\cal{G}}_i]B_{\theta(i)}du(i-d)\}\nonumber\\
&=&\mbox{E}\{2\sum_{i=d}^N\mbox{E}[\lambda_i'B_{\theta(i)}+u(i-d)'R_{\theta(i)}|{\cal{G}}_{i-d}]du(i-d)\}.\label{f16}
\end{eqnarray}
It concludes from (\ref{f16}) that the necessary condition for the minimum can be given as follows
\begin{eqnarray}
\mbox{E}[\lambda_i'B_{\theta(i)}+u(i-d)'R_{\theta(i)}|{\cal{G}}_{i-d}]=0,i=d,\cdots,N.\nonumber
\end{eqnarray}
(\ref{f3}) is shown. This completes the proof of Lemma 1.


\section{Proof of Proposition 1}

\emph{Proof}. In view of (\ref{a7})-(\ref{a9}), we have
\begin{eqnarray}
A_{l_k}'(\alpha_{l_k,l_k}^{d-1})'
&=&A_{l_k}'\{\Lambda_{l_k,l_{k+d-1}}[F_{l_{k+d-1},l_{k+1}}'(P_{l_{k+d-1}}-P_{l_{k+d-1}}^0)B_{l_{k+d-1}}]\nonumber\\
&&-\sum_{s=0}^{d-2}\Lambda_{l_k,l_{k+s}}[F_{l_{k+s},l_{k+1}}'(T_{l_{k+s}}^0)'W_{l_{k+s}}^{-1}T_{l_{k+s}}^{d-1-s}]\}.\label{fff5}
\end{eqnarray}
Recalling the definition of $T_{l_{k-1}}^0$, we get from (\ref{fff5}) that $\mbox{E}[A_{l_k}'(\alpha_{l_k,l_k}^{d-1})'|{\cal{G}}_{k-1}]=(T_{l_{k-1}}^0)'$.
(\ref{c1}) is shown. 
Similarly, we can show that (\ref{c2}), (\ref{c3}), and (\ref{c4})  are satisfied.


\section{Proof of Theorem 1}

\emph{Proof}. (i) Necessary: Assume that Problem 1 has a unique solution. We will prove that $W_{l_{k-d}}$ is invertible and $u(k-d)$ satisfies (\ref{a13}) for all $k=N,\cdots,d, l_{k-d}\in \Theta$ by the induction method. Define
\begin{eqnarray}
J(k)\stackrel{\triangle}{=} \mbox{E}\{\sum_{i=k}^N(x(i)'Q_{\theta(i)}x(i)+u(i-d)'R_{\theta(i)}u(i-d))+x(N+1)'P_{\theta(N+1)}x(N+1)|{\cal{G}}_{k-d}\},\label{f25}
\end{eqnarray}
for $k=N,\cdots,d$. For $k=N$, (\ref{f25}) becomes
\begin{eqnarray}
J(N){=} \mbox{E}\{x(N)'Q_{\theta(N)}x(N)+u(N-d)'R_{\theta(N)}u(N-d)+x(N+1)'P_{\theta(N+1)}x(N+1)|{\cal{G}}_{N-d}\}.\label{f26}
\end{eqnarray}
Based on (\ref{fff1}), we deduce that $J(N)$ can be formulated as a quadratic function of $x(N)$ and $u(N-d)$. The uniqueness of the optimal controller $u(N-d)$ indicates that the quadratic term of $u(N-d)$ is positive for any nonzero $u(N-d)$. Let $x(N)=0$ and substitute (\ref{fff1}) into (\ref{f26}), we have
\begin{eqnarray}
J(N)&=&\mbox{E}\{u(N-d)'(R_{\theta(N)}+B_{\theta(N)}'P_{\theta(N+1)}B_{\theta(N)})u(N-d)|{\cal{G}}_{N-d}\}\nonumber\\
&=&u(N-d)'W_{l_{N-d}}u(N-d)>0,l_{N-d}\in \Theta.\label{f27}
\end{eqnarray}
It can be concluded that $W_{l_{N-d}}>0$.

In what follows, the optimal controller $u(N-d)$ is to be calculated. Applying (\ref{fff1}), (\ref{f3}) and (\ref{f4}), we have
\begin{eqnarray}
0&=&\mbox{E}[B_{\theta(N)}'\lambda_N+R_{\theta(N)}u(N-d)|{\cal{G}}_{N-d}]\nonumber\\
&=&T_{l_{N-d}}^0x(N-d+1)+\sum_{j=1}^{d-1}T_{l_{N-d}}^ju(N-2d+j)+W_{l_{N-d}}u(N-d).\nonumber
\end{eqnarray}
It follows from the above equation that
\begin{eqnarray}
u(N-d)&=&W_{l_{N-d}}^{-1}T_{l_{N-d}}^0x(N-d+1)-\sum_{j=1}^{d-1}W_{l_{N-d}}^{-1}T_{l_{N-d}}^ju(N-2d+j).\label{a15}
\end{eqnarray}

In the following, we will show that $\lambda_{N-1}$ is with the form as (\ref{a14}). In view of (\ref{fff1}), (\ref{f5}), and (\ref{a15}), one yields
\begin{eqnarray}
\lambda_{N-1}
&=&P_{l_{N-1}}x(N)-(\alpha_{l_{N-1},l_{N-1}}^1)'W_{l_{N-d}}^{-1}[T_{l_{N-d}}^0x(N-d+1)-\sum_{j=1}^{d-1}T_{l_{N-d}}^ju(N-2d+j)].\nonumber
\end{eqnarray}
In view of the d-step backward formula (\ref{fff3}) and note that $P_{l_{N-1}}^0=0, (\alpha_{l_{N-1},l_{N-s}}^{d-s})'=0 (s=1,2,\cdots,d-2)$, one gets that
\begin{eqnarray}
\lambda_{N-1}&=&P_{l_{N-1}}x(N)-(\alpha_{l_{N-1},l_{N-1}}^1)'W_{l_{N-d}}^{-1}\mbox{E}\{\alpha_{l_{N-1},l_{N-1}}^1x(N)|{\cal{G}}_{N-d}\}\nonumber\\
&=&(P_{l_{N-1}}-P_{l_{N-1}}^0)x(N)-\sum_{s=1}^{d-1}(\alpha_{l_{N-1},l_{N-s}}^{d-s})'W_{l_{N-s-1}}^{-1}\mbox{E}\{\alpha_{l_{N-1},l_{N-s}}^{d-s}x(N)|{\cal{G}}_{N-s-1}\}.\nonumber
\end{eqnarray}
Thus we obtained (\ref{a14}) for $k=N$.

To proceed the induction proof, we take any $n$ with $1\leq n\leq N$, and assume that $W_{l_{k-d}}(k-d)$ is invertible and that the optimal controller $u(k-d)$ and the optimal costate$\lambda_{k-1}$ are as (\ref{a13}) and (\ref{a14}) for all $k\geq n+1$. In the next, it needs to show that these conditions will be satisfied for $k=n$. Let $x(n)=0$, we will check the quadratic term of $u(n-d)$ in $J(n)$. In view of (\ref{fff1}), (\ref{f3}) and (\ref{f5}) for $k\geq n+1$, we have
\begin{eqnarray}
&&\mbox{E}\{x(k)'\lambda_{k-1}-x(k+1)'\lambda_k|{\cal{G}}_{n-d+1}\}\}\nonumber\\
&=&\mbox{E}\{x(k)'Q_{\theta(k)}x(k)+u(k-d)'R_{\theta(k)}u(k-d)|{\cal{G}}_{n-d+1}\}.\nonumber
\end{eqnarray}
Adding from $k=n+1$ to $k=N$ on both sides of the above equation, we get
\begin{eqnarray}
&&\mbox{E}\{x(n+1)'\lambda_n-x(N+1)'\lambda_N|{\cal{G}}_{n-d+1}\}\nonumber\\
&=&\sum_{k=n+1}^N\mbox{E}\{x(k)'\lambda_{k-1}-x(k+1)'\lambda_k|{\cal{G}}_{n-d+1}\}\nonumber\\
&=&\sum_{k=n+1}^N\mbox{E}\{x(k)'Q_{\theta(k)}x(k)+u(k-d)'R_{\theta(k)}u(k-d)|{\cal{G}}_{n-d+1}\}.\label{f31}
\end{eqnarray}
It follows from (\ref{f31}) that
\begin{eqnarray}
J(n)
&=&\mbox{E}\{u(n-d)'R_{\theta(n)}u(n-d)+u(n-d)'B_{\theta(n)}'\lambda_n|{\cal{G}}_{n-d}\}.\label{f32}
\end{eqnarray}
Note that
\begin{eqnarray}
\lambda_{n}&=&(P_{l_{n}}-P_{l_{n}}^0)x(n+1)-\sum_{s=1}^{d-1}(\alpha_{l_{n},l_{n+1-s}}^{d-s})'W_{l_{n-s}}^{-1}\mbox{E}\{\alpha_{l_{n},l_{n+1-s}}^{d-s}x(n+1)|{\cal{G}}_{n-s}\}.\label{f33}
\end{eqnarray}
Substituting (\ref{f33}) into (\ref{f32}) and employing (\ref{c1})-(\ref{c4}), we get
\begin{eqnarray}
J(n)
&=&u(n-d)'W_{l_{n-d}}u(n-d).\label{a19}
\end{eqnarray}
It is concluded from the uniqueness of the optimal controller that $J(n)$ must be positive for any $u(n-d)\neq 0$. So we have $W_{l_{n-d}}(n-d)>0, l_{n-d}\in \Theta$.

To derive the optimal controller $u(n-d)$, plugging (\ref{f33}) in (\ref{f3}) yields
\begin{eqnarray}
0&=&\mbox{E}\{B_{\theta(n)}'\lambda_n+R_{\theta(n)}u(n-d)|{\cal{G}}_{n-d}\}\nonumber\\
&=&T_{l_{n-d}}^0x(n-d+1)+\sum_{j=1}^{d-1}T_{l_{n-d}}^ju(n-2d+j)+W_{l_{n-d}}u(n-d).\nonumber
\end{eqnarray}
Using the above equation, we get
\begin{eqnarray}
u(n-d)&=&-W_{l_{n-d}}^{-1}T_{l_{n-d}}^0x(n-d+1)-\sum_{j=1}^{d-1}W_{l_{n-d}}^{-1}T_{l_{n-d}}^ju(n-2d+j).\label{a23}
\end{eqnarray}
Now, we proceed to derive that $\lambda_{n-1}$ is of the form as (\ref{a14}). In terms of (\ref{f5}), (\ref{f33}) and (\ref{a23})and bearing in mind (\ref{c1})-(\ref{c4}) and (\ref{fff3}), we get
\begin{eqnarray}
\lambda_{n-1}&=&\mbox{E}\{Q_{\theta(n)}x(n)+A_{\theta(n)}'\lambda_n|{\cal{G}}_{n-1}\}\nonumber\\
&=&\mbox{E}\{Q_{l_n}x(n)+A_{l_n}'[(P_{l_n}-P_{l_n}^0)x(n+1)\nonumber\\
&&-\sum_{s=1}^{d-1}(\alpha_{l_n,l_{n+1-s}}^{d-s})'W_{l_{n-s}}^{-1}\mbox{E}(\alpha_{l_n,l_{n+1-s}}^{d-s}x(n+1)|{\cal{G}}_{n-s})]|{\cal{G}}_{n-1}\}\nonumber\\
&=&(P_{l_{n-1}}-P_{l_{n-1}}^0)x(n)-\sum_{s=1}^{d-1}(\alpha_{l_{n-1},l_{n-s}}^{d-s})'W_{l_{n-s-1}}^{-1}\mbox{E}\{\alpha_{l_{n-1},l_{n-s}}^{d-s}x(n)|{\cal{G}}_{n-s-1}\}.\label{a24}
\end{eqnarray}

(\ref{a14}) is proved. The proof of necessity is finished.

(ii) Sufficiency-The proof of sufficiency is similar to that of Theorem 1 in \cite{ZLXF15}.
This completes the proof Theorem 1.


\end{document}